\documentclass[12pt,a4paper,oneside,reqno]{article}

\usepackage[utf8]{inputenc}
\usepackage[english]{babel}
\usepackage[left=1in,right=1in,top=1in,bottom=1in]{geometry}
\usepackage{amsfonts,amssymb,amsmath,amsthm}
\usepackage{xcolor}
\usepackage[normalem]{ulem}
\usepackage{cite}
\usepackage{hyperref}
\hypersetup{colorlinks=true,linkcolor=blue,citecolor=blue,urlcolor=blue}
\usepackage{fancyhdr}
\newtheorem{theorem}{Theorem}
\newtheorem{lemma}[theorem]{Lemma}

\newcommand{\doi}[1]{\href{https://doi.org/#1}{doi:#1}}
\newcommand{\arxiv}[1]{\href{https://arxiv.org/abs/#1}{arXiv:#1}}

\begin{document}

\noindent\fbox{\parbox{\dimexpr\textwidth-2\fboxsep-2\fboxrule\relax}{\small
\textbf{Corrected version verified in Lean.}  This is a corrected version of the paper
I.~Kaygorodov, E.~Okhapkina, \emph{$\delta$-derivations of semisimple finite-dimensional
structurable algebras}, J. Algebra Appl. 13 (2014), no.~4, 1350130,
\doi{10.1142/S0219498813501302}.}}

\bigskip

\begin{center}
{\Large $\delta$-derivations of semisimple\\[2pt]  finite-dimensional structurable algebras}

\bigskip

{\bf Ivan Kaygorodov and Elizaveta Okhapkina}

\end{center}

\medskip

\noindent{\small\textbf{Abstract.} In this paper, we show the absence of nontrivial
$\delta$-derivations of semisimple finite-dimensional structurable algebras over an
algebraically closed field of characteristic not equal to $2,3,5$.}

\medskip


\bigskip


\bigskip

\setcounter{tocdepth}{1}
\tableofcontents

\section{Introduction}

The concept of $\delta$-derivation first appeared in papers by Filippov~\cite{3,4,5}, as a
generalization of ordinary derivations. Recall that for a fixed $\delta$ of the main field $F$,
under a $\delta$-derivation of an algebra $A$ we understand a linear mapping $\varphi$ which
satisfies the condition
\begin{equation}\label{e1}
\varphi(xy)=\delta(\varphi(x)y+x\varphi(y))
\end{equation}
for arbitrary elements $x,y\in A$. In paper~\cite{4}, Filippov proved that every prime Lie
algebra has no nonzero $\delta$-derivations if $\delta\neq-1,0,\frac12,1$. Subsequently, in
papers~\cite{6,12,17} the description of $\delta$-derivations of semisimple finite-dimensional
Jordan algebras and $\delta$-(super)derivations of semisimple finite-dimensional Jordan
superalgebras was given. Later, Shestakov generalized some results of~\cite{6} and described
ternary derivations of simple finite-dimensional Jordan algebras~\cite{21}. $\delta$-derivations
of Lie superalgebras were considered by Kaygorodov~\cite{7,8} and Zusmanovich~\cite{26}. The
analogue of $\delta$-derivation for $n$-ary algebras was investigated in~\cite{13,14}. We note
that some generalizations of $\delta$-derivations were proposed in
papers~\cite{9,10,15,16,24}, and $\delta$-derivations were used for studying prime Lie algebras
by Filippov~\cite{3}, for the construction of post-Lie structures by Burde and Dekimpe~\cite{2},
for studying speciality of Jordan superalgebras by Kaygorodov and Zhelyabin~\cite{17}, in the
analysis of the derivations of current Lie algebras by Zusmanovich~\cite{27}, and for studying
homotopes of Novikov algebras by Filippov and Sereda~\cite{20}.

\section{Examples and definitions}

The class of structurable algebras was introduced in 1978 by Allison~\cite{1}. It contains such
objects as tensor products of composition algebras, the $56$-dimensional Freudenthal module over
$E_7$ with a natural binary operation, and a $35$-dimensional algebra $T(C)$. Also, for these
algebras a natural generalization of the Kantor--Koecher--Tits construction holds.

Structurable algebras are algebras with unit, which we are going to designate as $e$, and
involution $\bar{\ }$, which satisfy the identity
\begin{equation}\label{e2}
[T_z,V_{x,y}]=V_{T_zx,y}-V_{x,T_{\bar z}y},
\end{equation}
where $T_z,V_{x,y}\in\mathrm{End}(A)$,
\[
V_{x,y}(z)=(x\bar y)z+(z\bar y)x-(z\bar x)y,\qquad T_z=V_{z,1}
\]
for $x,y,z\in A$.

Let $(A,\bar{\ })$ be an algebra with involution. We can consider the vector space $A$ as a direct
sum $A=H\oplus S$, where $H=\{a\in A\mid \bar a=a\}$ is the set of symmetric elements and
$S=\{a\in A\mid \bar a=-a\}$ is the set of antisymmetric elements.

Let us give some examples of structurable algebras (see~\cite{22,23}).

\smallskip

\noindent(S1) Let $(A,\bar{\ })$ be a unital associative algebra with involution. Then
$(A,\bar{\ })$ is a structurable algebra.

\smallskip

\noindent(S2) Let $J$ be a Jordan algebra and $\bar{\ }=\mathrm{id}$. For $(J,\bar{\ })$ the
identity~\eqref{e2} turns into the identity
\[
[L_c,V_{a,b}]=V_{ca,b}-V_{a,cb},
\]
which holds in an arbitrary Jordan algebra. Therefore $(J,\bar{\ })$ is a structurable algebra.

\smallskip

\noindent(S3) Let $(E,\bar{\ })$ be an associative algebra with involution and $W$ a left
$E$-module relative to the operation $\circ\colon E\times W\to W$. Let $W$ possess a Hermitian
form $h$, that is, a bilinear mapping $h\colon W\times W\to E$ such that
\[
\overline{h(w_1,w_2)}=h(w_2,w_1),\qquad h(e_1\circ w_1,w_2)=e_1h(w_1,w_2)
\]
for any $w_1,w_2\in W$ and $e_1\in E$. The operations of multiplication and involution on the
space $A=E\oplus W$ are defined as follows:
\[
(e_1,w_1)\bullet(e_2,w_2)=(e_1e_2+h(w_2,w_1),\ e_2\circ w_1+\bar e_1\circ w_2),\qquad
\overline{(e,w)}=(\bar e,w).
\]
The algebra $(A,\bar{\ })$ is a structurable algebra. It is a generalization of the Jordan algebra
of a bilinear form and is called the algebra of a Hermitian form. The algebra of a Hermitian form
is simple if the Hermitian form $h$ on the left $E$-module is nondegenerate and $(E,\bar{\ })$ is
a central simple algebra with involution.

\smallskip

\noindent(S4) Let $(G_1,\bar{\ })$ and $(G_2,\bar{\ })$ be composition algebras with standard
involution (for more about composition algebras see, for example,~\cite{25}). Let us define on
the algebra $A=G_1\otimes G_2$ the involution $\overline{x_1\otimes x_2}=\bar x_1\otimes\bar x_2$.
Then $(A,\bar{\ })$ is a structurable algebra.

\smallskip

\noindent(S5) Let $J,J'$ be vector spaces over $F$ with trilinear forms $N$ and $N'$ which are
related by a nondegenerate bilinear form $T\colon J\times J'\to F$. For $j,k\in J$,
$j',k'\in J'$, we define $j\times k\in J'$ and $j'\times k'\in J$ by the conditions
\[
T(l,j\times k)=N(j,k,l),\qquad T(j'\times k',l')=N'(j',k',l')
\]
for any $l\in J$, $l'\in J'$. On $J$ and $J'$ we define the operation $\#$ as
\[
j^\#=\tfrac12\,j\times j,\qquad j'^\#=\tfrac12\,j'\times j'\qquad\text{for } j\in J,\ j'\in J'.
\]
The triple $(T,N,N')$ is called an admissible triple defined on the pair of spaces $(J,J')$ if $N$
and $N'$ are nontrivial and satisfy the conditions
\[
(j^\#)^\#=\tfrac16N(j,j,j)j,\qquad (j'^\#)^\#=\tfrac16N'(j',j',j')j'\qquad
\text{for } j\in J,\ j'\in J'.
\]
Let $(T,N,N')$ be an admissible triple on a pair of spaces $(J,J')$. On the space of matrices
\[
A=\left\{\begin{pmatrix}\alpha&j\\ j'&\beta\end{pmatrix}:\ \alpha,\beta\in F,\ j\in J,\ j'\in J'\right\}
\]
the multiplication
\[
\begin{pmatrix}\alpha&j\\ j'&\beta\end{pmatrix}\begin{pmatrix}\gamma&k\\ k'&\delta\end{pmatrix}
=\begin{pmatrix}\alpha\gamma+T(j,k')&\alpha k+\delta j+j'\times k'\\
\gamma j'+\beta k'+j\times k&T(k,j')+\beta\delta\end{pmatrix}
\]
and the involution
\[
\overline{\begin{pmatrix}\alpha&j\\ j'&\beta\end{pmatrix}}=\begin{pmatrix}\beta&j\\ j'&\alpha\end{pmatrix}
\]
are defined. The obtained structurable algebra $(A,\bar{\ })$ is called the algebra constructed
from an admissible triple $(T,N,N')$ on a pair of spaces $(J,J')$.

\smallskip

\noindent(S6) Let $C$ be a Cayley--Dickson algebra with involution (about Cayley--Dickson algebras
see~\cite{25}). The set $S$ of antisymmetric elements is a seven-dimensional simple non-Lie Malcev
algebra relative to the commutator $[\cdot,\cdot]$ in $C$ (see~\cite{18}); it is $S^{(-)}$. On
$S^{(-)}$ there is a symmetric nondegenerate invariant bilinear form $(\cdot,\cdot)$ such that for
$x,y\in S$ the condition
\[
[[x,y],y]=(y,y)x-(x,y)y
\]
is satisfied. Let $M$ be the subspace of the tensor product $S\otimes S$ which is generated by the
set $\{s\otimes r-r\otimes s: s,r\in S\}$. We let $H=S\otimes S/M$, and on the direct sum of spaces
$H\oplus S$ the commutative $\odot$ and anticommutative $[\cdot,\cdot]$ operations are defined by
the following conditions:
\begin{align*}
[s_1,s_2]&=[s_1,s_2],\\
[s,s_1\otimes s_2]&=[s,s_1]\otimes s_2+s_1\otimes[s,s_2],\\
[s_1\otimes s_2,s_3\otimes s_4]&=(s_1,s_3)[s_2,s_4]+(s_1,s_4)[s_2,s_3]+(s_2,s_3)[s_1,s_4]
+(s_2,s_4)[s_1,s_3],\\
s_1\odot s_2&=s_1\otimes s_2,\\
s\odot(s_1\otimes s_2)&=\tfrac12(s_1,s_2)s+\tfrac14(s,s_1)s_2+\tfrac14(s,s_2)s_1,\\
(s_1\otimes s_2)\odot(s_3\otimes s_4)&=\tfrac14[s_1,s_3]\otimes[s_2,s_4]
+\tfrac14[s_1,s_4]\otimes[s_2,s_3]\\
&\quad+\tfrac12(s_1,s_2)s_3\otimes s_4+\tfrac12(s_3,s_4)s_1\otimes s_2
\end{align*}
for $s,s_1,s_2,s_3,s_4\in S$, where on the right-hand side under $[\cdot,\cdot]$ we understand the
commutator in the algebra $C$. The operations of multiplication and involution in $H\oplus S$ are
defined in the following way:
\[
xy=x\odot y+\tfrac12[x,y],\qquad \overline{h+s}=h-s,\qquad\text{where } x,y\in H\oplus S,\
h\in H,\ s\in S.
\]
The algebra with involution so obtained, which is built by means of the Cayley--Dickson algebra
$C$, is designated by $T(C)$.

\smallskip

Let us note that, according to the results of~\cite{22,23}, the algebras of types (S1)--(S6)
exhaust all simple finite-dimensional structurable algebras over algebraically closed fields of
characteristic not equal to $2,3,5$.

As was noted above, for a fixed element $\delta$ of the main field, under a $\delta$-derivation of
an algebra $A$ we understand a linear mapping $\varphi\colon A\to A$ which for arbitrary $x,y\in A$
satisfies the condition
\[
\varphi(xy)=\delta(\varphi(x)y+x\varphi(y)).
\]
The centroid $\Gamma(A)$ of an algebra $A$ is the set of linear mappings $\chi\colon A\to A$ with
the condition $\chi(ab)=\chi(a)b=a\chi(b)$. It is clear that any element of the centroid of an
algebra is a $\frac12$-derivation. Any endomorphism $\varphi$ of an algebra $A$ such that
$\varphi(A^2)=0$ is a $0$-derivation. Any derivation of an algebra $A$ is a $1$-derivation. A
nonzero $\delta$-derivation $\varphi$ is a nontrivial $\delta$-derivation if $\delta\neq0,1$ and
$\varphi\notin\Gamma(A)$.

\section{Main lemmas}

Here and below, all algebras are considered over an algebraically closed field $F$ of
characteristic not equal to $2,3,5$, although some of the results of the following lemmas remain
correct under weaker conditions on the field.

Making a trivial direct calculation and using~\cite[Theorem~2.1]{6}, we obtain that if $\varphi$
is a nontrivial $\delta$-derivation of a structurable algebra $A$, then $\delta=\frac12$ and for
any element $x\in A$ we have $\varphi(x)=ax=xa$, where $a\in A$ is a certain fixed element.

\begin{lemma}\label{l1}
Let $A$ be the algebra of a nondegenerate Hermitian form on a left $E$-module, where
$(E,\bar{\ })$ is a simple central associative algebra with involution. Then every
$\frac12$-derivation of $A$ is trivial.
\end{lemma}

\begin{proof}
Let $\varphi$ be a $\frac12$-derivation of $A$ and let $\varphi(e)=(c,u)$, where $c\in E$, $u\in W$.
Then $\varphi(X)=X\bullet\varphi(e)=\varphi(e)\bullet X$ for all $X\in A$. Comparing
\begin{center}$\varphi((x,0))=(c,u)\bullet(x,0)=(cx,x\circ u)$
\end{center} with
$\varphi((x,0))=(x,0)\bullet(c,u)=(xc,\bar x\circ u)$, we get
\begin{equation}\label{e3}
cx=xc,\qquad x\circ u=\bar x\circ u\qquad\text{for all } x\in E.
\end{equation}
Comparing $\varphi((0,w))=(c,u)\bullet(0,w)=(h(w,u),\bar c\circ w)$ with\begin{center}
$\varphi((0,w))=(0,w)\bullet(c,u)=(h(u,w),c\circ w)$,\end{center} we get $h(w,u)=h(u,w)$ for all $w\in W$.

Let $K=\{z\in E\mid z\circ u=0\}$. Clearly, $K$ is a left ideal of $E$. If $z\in K$ and $y\in E$,
then by~\eqref{e3}
\[
(zy)\circ u=\overline{zy}\circ u=\bar y\circ(\bar z\circ u)=\bar y\circ(z\circ u)=0,
\]
so $K$ is a two-sided ideal, and $\bar z\circ u=z\circ u=0$ shows that $\bar K\subseteq K$. Since
$(E,\bar{\ })$ is simple, either $K=E$ or $K=0$.

If $K=E$, then $u=e\circ u=0$ and $\varphi(X)=(c,0)\bullet X$, where $c$ lies in the center of $E$
by~\eqref{e3} (hence so does $\bar c$). For $X=(x,v)$, $Y=(y,w)$ we have
\begin{center}$h(w,\bar c\circ v)=\overline{\bar ch(v,w)}=h(w,v)c=ch(w,v)$,\end{center} and therefore
\begin{align*}
\varphi(X\bullet Y)&=(c,0)\bullet(xy+h(w,v),\ y\circ v+\bar x\circ w)\\
&=(cxy+ch(w,v),\ \bar c\circ(y\circ v)+\bar c\circ(\bar x\circ w))\\
&=(cxy+h(w,\bar c\circ v),\ y\circ(\bar c\circ v)+\overline{cx}\circ w)
=(cx,\bar c\circ v)\bullet(y,w)=\varphi(X)\bullet Y.
\end{align*}
Since $\varphi$ is a $\frac12$-derivation, this also gives $\varphi(X\bullet Y)=X\bullet\varphi(Y)$,
i.e., $\varphi\in\Gamma(A)$.

If $K=0$, then $x=\bar x$ for all $x\in E$ by~\eqref{e3}; hence $E$ is commutative
\begin{center}($xy=\overline{xy}=\bar y\bar x=yx$)\end{center} and $h$ is symmetric. Applying the definition of a
$\frac12$-derivation to $X=(0,v)$, $Y=(0,w)$ and comparing the components in $W$, we get
\begin{equation}\label{e4}
2h(w,v)\circ u=h(u,v)\circ w+h(u,w)\circ v\qquad\text{for all } v,w\in W.
\end{equation}
For $w=u$ this gives $h(u,v)\circ u=t\circ v$, where $t=h(u,u)$. If $t=0$, then $h(u,v)\in K=0$ for
all $v\in W$, and $u=0$ by nondegeneracy of $h$, which contradicts $K=0$. Thus $Et$ is a nonzero
ideal of $E$ stable under the (identity) involution, so $Et=E$ and $t$ is invertible. Then
$v=t^{-1}h(u,v)\circ u$, whence $h(u,v)\circ w=t^{-1}h(u,v)h(u,w)\circ u=h(u,w)\circ v$, and
\eqref{e4} becomes $h(w,v)\circ u=h(v,u)\circ w$. Now for $X=(x,v)$, $Y=(y,w)$ we obtain
{\small \begin{align*}
\varphi(X\bullet Y)&=(c,u)\bullet(xy+h(w,v),\ y\circ v+x\circ w)\\
&=(cxy+ch(w,v)+yh(v,u)+xh(w,u),\ xy\circ u+h(w,v)\circ u+cy\circ v+cx\circ w)\\
&=(cxy+h(v,u)y+h(w,x\circ u)+h(w,c\circ v),\ y\circ(x\circ u+c\circ v)+(cx+h(v,u))\circ w)\\
&=(cx+h(v,u),\ x\circ u+c\circ v)\bullet(y,w)=\varphi(X)\bullet Y,
\end{align*}}
and, as above, $\varphi\in\Gamma(A)$.
The lemma is proved.
\end{proof}

\begin{lemma}\label{l2}
The algebra of the tensor product of composition algebras has no nontrivial
$\frac12$-derivations.
\end{lemma}

\begin{proof}
Let $A,B$ be composition algebras (about composition algebras see, for example,~\cite{25}) with
units $1_A,1_B$, and let $e$ be the unit of the algebra $A\otimes B$. If $\varphi$ is a
$\frac12$-derivation of the algebra $A\otimes B$, then $\varphi(e)=\sum_ia_i\otimes b_i$,
$a_i\in A$, $b_i\in B$, where the elements $b_i$ are linearly independent. Since $(a\otimes1_B)\varphi(e)=\varphi(e)(a\otimes1_B)$ for $a\in A$, the
elements $a_i$ are in the commutative center of $A$.

It is easy to see that
\[
\varphi(ac\otimes1_B)=\tfrac12\varphi(a\otimes1_B)(c\otimes1_B)+\tfrac12(a\otimes1_B)\varphi(c\otimes1_B).
\]
From the previous relation we get
\begin{multline*}
(ac\otimes1_B)\varphi(1_A\otimes1_B)=\varphi(ac\otimes1_B)=\varphi((a\otimes1_B)(c\otimes1_B))\\
=\tfrac12((a\otimes1_B)\varphi(1_A\otimes1_B))(c\otimes1_B)
+\tfrac12(a\otimes1_B)((c\otimes1_B)\varphi(1_A\otimes1_B)),
\end{multline*}
that is,
\[
\sum_i(ac)a_i\otimes b_i=\tfrac12\sum_i\bigl((aa_i)c\otimes b_i+a(ca_i)\otimes b_i\bigr).
\]
The last relation implies
\[
\sum_i\bigl(2(ac)a_i-(aa_i)c-a(ca_i)\bigr)\otimes b_i=0\quad\text{and}\quad
2(ac)a_i-(aa_i)c-a(ca_i)=0.
\]
Since $a_i$ is in the commutative center of $A$, the last relation
can be rewritten as $2(a,c,a_i)-(a,a_i,c)=0$. As composition algebras are alternative,
$(a,a_i,c)=-(a,c,a_i)$, whence $3(a,c,a_i)=0$ and $(a,c,a_i)=0$; therefore all associators
containing $a_i$ vanish. Thus, $a_i$ is in the commutative--associative center $Z(A)$, i.e.,
$\varphi(e)\in Z(A)\otimes B$.

As the algebra $A\otimes B$ is isomorphic to the algebra $B\otimes A$, analogously we get
that
$\varphi(e)\in A\otimes Z(B)$. Hence $\varphi(e)\in(Z(A)\otimes B)\cap(A\otimes Z(B))=Z(A)\otimes
Z(B)$, i.e., we may assume that $a_i\in Z(A)$ and $b_i\in Z(B)$.

Now we are going to show that the mapping $\varphi\colon A\otimes B\to A\otimes B$ which is
defined by the rule
\[
\varphi(x\otimes y)=(x\otimes y)\Bigl(\sum_ia_i\otimes b_i\Bigr),\qquad\text{where } a_i\in Z(A),\
b_i\in Z(B),
\]
is an element of the centroid of the algebra $A\otimes B$, that is, a trivial
$\frac12$-derivation. It is easy to see that
\begin{multline*}
\varphi((a\otimes b)(c\otimes d))=\varphi(ac\otimes bd)=(ac\otimes bd)\Bigl(\sum_ia_i\otimes b_i\Bigr)\\
=\sum_i(ac)a_i\otimes(bd)b_i=\sum_i(aa_i\otimes bb_i)(c\otimes d)=\varphi(a\otimes b)(c\otimes d).
\end{multline*}
Analogously, we show that $\varphi((a\otimes b)(c\otimes d))=(a\otimes b)\varphi(c\otimes d)$. The
lemma is proved.
\end{proof}

\begin{lemma}\label{l3}
The algebra of an admissible triple $(T,N,N')$ has no nontrivial $\frac12$-derivations.
\end{lemma}

\begin{proof}
Let $\varphi(e)=\begin{pmatrix}\alpha&j\\ j'&\beta\end{pmatrix}$ and
$x=\begin{pmatrix}\gamma&k\\ k'&\delta\end{pmatrix}$, where $\alpha,\beta,\gamma,\delta\in F$,
$k,j\in J$, $k',j'\in J'$. We have $x\varphi(e)=\varphi(e)x$, so
\[
\varphi(e)x=\begin{pmatrix}\alpha&j\\ j'&\beta\end{pmatrix}\begin{pmatrix}\gamma&k\\ k'&\delta\end{pmatrix}
=\begin{pmatrix}\alpha\gamma+T(j,k')&\alpha k+\delta j+j'\times k'\\
\gamma j'+\beta k'+j\times k&T(k,j')+\beta\delta\end{pmatrix}
\]
and
\[
x\varphi(e)=\begin{pmatrix}\gamma&k\\ k'&\delta\end{pmatrix}\begin{pmatrix}\alpha&j\\ j'&\beta\end{pmatrix}
=\begin{pmatrix}\gamma\alpha+T(k,j')&\gamma j+\beta k+k'\times j'\\
\alpha k'+\delta j'+k\times j&T(j,k')+\delta\beta\end{pmatrix}.
\]
Thus, it is necessary that
\[
\gamma j'+\beta k'+j\times k=\alpha k'+\delta j'+k\times j
\]
and
\[
\alpha k+\delta j+j'\times k'=\gamma j+\beta k+k'\times j'.
\]
For $\gamma=1$, $\delta=0$, $k=0$,
$k'=0$ the first condition gives $j'=0$; for $\gamma=0$, $\delta=1$, $k=0$, $k'=0$ the second
condition gives $j=0$. Since $N$ is nontrivial, $J\neq0$; taking $\gamma=\delta=0$, $k'=0$ and
$0\neq k\in J$ in the second condition, we get $\alpha k=\beta k$, i.e., $\alpha=\beta$. Hence,
$\varphi(e)=\alpha\begin{pmatrix}1&0\\0&1\end{pmatrix}$. The lemma is proved.
\end{proof}

\begin{lemma}\label{l4}
The algebra $T(C)$ has no nontrivial $\frac12$-derivations.
\end{lemma}

\begin{proof}
Let $e_1,\dots,e_8$ be the standard basis of the algebra $C$ (see, for example,~\cite{25}). Let us
recall that $S^{(-)}$ is a simple seven-dimensional non-Lie Malcev algebra relative to the
commutator $[\cdot,\cdot]$ in $C$ (about Malcev algebras see, for example,~\cite{18}), where $S$ is
the set of elements of the algebra $C$ that are antisymmetric relative to the involution. It is easy
to see that $e_i\in S$ if $i=2,3,\dots,8$.

Let $\varphi(e)=h^*+s^*$, where $h^*\in H$, $s^*\in S$, and $e$ is the unit of the algebra $T(C)$.
Then
\[
[\varphi(e),s]=0\quad\text{and}\quad[h^*,s]+[s^*,s]=0.
\]
Let us notice that $[s^*,s]=0$. As the algebra $S^{(-)}$ is simple, its annihilator is trivial. It
is easy to see that $s^*$ is in the annihilator of the algebra $S^{(-)}$, that is, $s^*=0$.
Consequently, $[h^*,s]=0$.

It is easy to see (see, for example,~\cite{19}) that in the Cayley--Dickson algebra $C$ a basis can
be chosen in such a way that the multiplication of the basis elements of the algebra is represented
by the following table:
\begin{align*}
e_2&=e_1e_2=e_5e_6=e_7e_8=e_3e_4, & e_3&=e_1e_3=e_7e_6=e_4e_2=e_8e_5,\\
e_4&=e_1e_4=e_2e_3=e_6e_8=e_7e_5, & e_5&=e_1e_5=e_6e_2=e_4e_7=e_3e_8,\\
e_6&=e_1e_6=e_2e_5=e_8e_4=e_3e_7, & e_7&=e_1e_7=e_5e_4=e_8e_2=e_6e_3,\\
e_8&=e_1e_8=e_2e_7=e_4e_6=e_5e_3. &&
\end{align*}
Let the element $h^*$ be expressed by the basis elements in the following way:
\[
h^*=\sum_{i,j=2}^8\alpha_{ij}e_i\otimes e_j.
\]
We can say that $[h^*,x]=0=[x,h^*]$. Let us analyze the product $[e_j,h^*]$. From the condition
\[
\Bigl[e_2,\sum_{i,j=2}^8\alpha_{ij}e_i\otimes e_j\Bigr]=0
\]
we get conditions on the $\alpha_{ij}$, and then the element of interest can be represented as
follows:
\begin{multline*}
h^*=\alpha_{22}e_2\otimes e_2+\alpha_{33}e_3\otimes e_3+\alpha_{33}e_4\otimes e_4
+\alpha_{35}e_3\otimes e_5+\alpha_{45}e_4\otimes e_5+\alpha_{55}e_5\otimes e_5\\
-\alpha_{45}e_3\otimes e_6+\alpha_{35}e_4\otimes e_6+\alpha_{55}e_6\otimes e_6
+\alpha_{37}e_3\otimes e_7+\alpha_{47}e_4\otimes e_7+\alpha_{57}e_5\otimes e_7
+\alpha_{67}e_6\otimes e_7\\
+\alpha_{77}e_7\otimes e_7-\alpha_{47}e_3\otimes e_8
+\alpha_{37}e_4\otimes e_8-\alpha_{67}e_5\otimes e_8+\alpha_{57}e_6\otimes e_8
+\alpha_{77}e_8\otimes e_8.
\end{multline*}
Performing the multiplication by $e_j$ ($j=3,4,5,6,7$), we get that all $\alpha_{ij}=0$ when
$i\neq j$. Combining these results with those obtained from the condition $[e_8,h^*]=0$, it is easy
to notice that $\alpha_{ii}=\alpha_{jj}$ for any $i,j$. Therefore,
$h^*=\sum_{i=2}^8\alpha e_i\otimes e_i$.

Let us show that $(e_i,e_i)=-4$ and $(e_i,e_j)=0$ when $i\neq j$:
\[
[[e_j,e_i],e_i]=(e_i,e_i)e_j-(e_j,e_i)e_i;
\]
\[
[e_j,e_i]e_i-e_i[e_j,e_i]=2[e_j,e_i]e_i=2(e_je_i-e_ie_j)e_i=4(e_je_i)e_i=-4e_j.
\]
Using this fact and the multiplication in the algebra $T(C)$, by direct substitution and
verification we can easily see that the element $e=-\frac1{16}\sum_{i=2}^8e_i\otimes e_i$ is the
unit of the algebra. Let $\beta=-16\alpha$. Then $h^*=\beta e$. The
lemma is proved.
\end{proof}

\section{Main results}

\begin{theorem}\label{t5}
A simple finite-dimensional structurable algebra over an algebraically closed field of
characteristic $p\neq2,3,5$ has no nontrivial $\delta$-derivations.
\end{theorem}

\begin{proof}
Following~\cite{23}, we can say that a simple finite-dimensional structurable algebra over an
algebraically closed field of characteristic $p\neq2,3,5$ is isomorphic to one of the algebras of
the types (S1)--(S6) listed in Section~2. Let us notice that nontrivial $\delta$-derivations are
possible only when $\delta=\frac12$. We have that
\begin{itemize}
\item[(1)] algebras of types (S1), (S2) do not admit nontrivial $\frac12$-derivations
(see~\cite{6,12}); for type (S1), if $a=\varphi(e)$, then
\begin{center}$\varphi(xy)=a(xy)=(ax)y=\varphi(x)y$ and $\varphi(xy)=(xy)a=x(ya)=x\varphi(y)$,  \end{center}so
$\varphi\in\Gamma(A)$;
\item[(2)] algebras of types (S3)--(S6) do not admit nontrivial $\frac12$-derivations (see
Lemmas~\ref{l1}--\ref{l4}).
\end{itemize}
Thus, the theorem is proved.
\end{proof}

\begin{theorem}\label{t6}
A semisimple finite-dimensional structurable algebra over an algebraically closed field of
characteristic $p\neq2,3,5$ has no nontrivial $\delta$-derivations.
\end{theorem}

\begin{proof}
According to~\cite{23}, if $A$ is a semisimple finite-dimensional structurable algebra over an
algebraically closed field of characteristic not equal to $2,3,5$, then $A=\bigoplus A_i$, where
$A_i$ is a simple structurable algebra. Let $e_k$ be the unit of the algebra $A_k$, and let
$\varphi$ be a $\delta$-derivation of $A$, $\delta\neq0,1$. If $x_i\in A_{i}$, then
$\varphi(x_i)=x_i^++x_i^-$, where $x_i^+\in A_i$, $x_i^-\notin A_i$ lies in the sum of the
other components. Let $e^i=\sum_{k\neq i}e_k$ and $\varphi(e^i)=e^{i+}+e^{i-}$, where
$e^{i+}\in A_i$, $e^{i-}\notin A_i$. Then
\[
0=\varphi(x_ie^i)=\delta(\varphi(x_i)e^i+x_i\varphi(e^i))
=\delta((x_i^++x_i^-)e^i+x_i(e^{i+}+e^{i-}))=\delta(x_i^-+x_ie^{i+}),
\]
whence $x_i^-=0$. Thus $\varphi$ is invariant on $A_i$. According to
Theorem~\ref{t5}, $A_i$ has no nontrivial $\delta$-derivations, so the restriction
$\varphi_i$ of $\varphi$ to $A_i$ belongs to $\Gamma(A_i)$. Since the multiplication in $A$ is
componentwise, $\varphi=\bigoplus\varphi_i\in\Gamma(A)$. Then the semisimple algebra $A$ has no
nontrivial $\delta$-derivations. The theorem is proved.
\end{proof}

Let us note that in~\cite{11} the definition of a generalized $\delta$-derivation is introduced. A
linear mapping $\chi$ is called a generalized $\delta$-derivation if it is related to a
$\delta$-derivation $\varphi$ by the condition
\[
\chi(xy)=\delta(\chi(x)y+x\varphi(y))=\delta(\varphi(x)y+x\chi(y)).
\]
According to the results of~\cite{11}, unital algebras have no generalized $\delta$-derivations
that are neither generalized derivations nor $\delta$-derivations. It was also shown that if a
generalized $\delta$-derivation is a generalized derivation, then it is the sum of a derivation and
an element of the centroid. Therefore, from Theorem~\ref{t6} the theorem below follows.

\begin{theorem}\label{t7}
Let $\chi$ be a generalized $\delta$-derivation of a semisimple finite-dimensional structurable
algebra $A$ over an algebraically closed field of characteristic $p\neq2,3,5$. Then
$\chi\in\mathrm{Der}(A)+\Gamma(A)$, where $\mathrm{Der}(A)$ is the space of derivations and
$\Gamma(A)$ is the centroid.
\end{theorem}

\section*{Appendix A. Differences from the published version}
\addcontentsline{toc}{section}{Appendix A. Differences from the published version}

\emph{Lemma~\ref{l1}: proof changed.}

\noindent\emph{Lemmas~\ref{l2}--\ref{l4}, Theorems~\ref{t5} and~\ref{t6}}: proofs corrected.

\end{document}